# Superefficiency from the Vantage Point of Computability

## Vladimir Vovk

*Abstract.* In 1952 Lucien Le Cam announced his celebrated result that, for regular univariate statistical models, sets of points of superefficiency have Lebesgue measure zero. After reviewing the turbulent history of early studies of superefficiency, I suggest using the notion of computability as a tool for exploring the phenomenon of superefficiency. It turns out that only computable parameter points can be points of superefficiency for computable estimators. This algorithmic version of Le Cam's result implies, in particular, that sets of points of superefficiency not only have Lebesgue measure zero but are even countable.

... if the true [parameter] value were known, a system of estimation could be devised which would give it with arbitrarily small variance; and such a system of estimation might happen to be adopted even if the true value were unknown.

Harold Hotelling, from a letter to
Ronald A. Fisher, 1930

## 1. INTRODUCTION

At the beginning of his recent paper [45] Stephen Stigler presents Hodges's famous example of a superefficient estimator as a nasty, ugly little fact that killed Fisher's beautiful theory of efficiency of maximum likelihood. Extending and permuting Wolfowitz's [52] classification, we call the three main lines of defense against the little fact "exclusion of the evil" (eliminating the superefficient estimators from competition), "deprecation of the evil" (showing that superefficiency can only happen on a small set of parameter points) and "collective responsibility" (refusing to accept a parameter point as a point of superefficiency, or even simply efficiency, unless its neighbors are points of efficiency). They will be reviewed in Section 2. Our review will be rather selective and will end around 1970—by that time the theory of superefficiency for regular parametric models had been essentially completed.

*Vladimir Vovk is Professor, Department of Computer Science, Royal Holloway, University of London, Egham, Surrey TW20 0EX, UK e-mail: vovk@cs.rhul.ac.uk.*

The rest of this paper concentrates on the second line of defense, with a minimal, and very natural, admixture of the first line: we will restrict our attention to the computable estimators. On the other hand, we will never assume asymptotic normality of our estimators, although our definition of asymptotic efficiency is motivated by comparison with the asymptotically normal case. The result that superefficiency can occur only at computable parameter points is established in Section 3 as Theorem 1. Surprisingly, the regularity conditions required for this result are relatively simple and easy to check; this is discussed in Section 4.

The notion of computability for real numbers and functions of real numbers will be defined and discussed in Appendix A. The proof of Theorem 1 will be very brief in the part concerning computability, and the details will be provided in Appendix A.

Another appendix, Appendix B, contains a direct proof of the countability of sets of superefficiency not using the notion of computability.

The absence of superefficiency for computable continuous estimators at noncomputable points was first established in [49] in the framework of Kolmogorov's algorithmic theory of randomness (see, e.g., [24]). A serious disadvantage of the algorithmic theory of randomness is its unfamiliarity to most statisticians. Another disadvantage is that typical results proved in the framework of the algorithmic theory of randomness contain unspecified constants, which mask important details. This paper allows noncontinuous estimators and avoids using the algorithmic theory of randomness.

There is more than one connection between superefficiency and computing. This paper applies the notion of computability to study superefficiency. In the opposite direction, Barron and Hengartner [4] use the notion of superefficiency to study the important computational problem of data compression.

## 2. FISHER'S PROGRAM AND SUPEREFFICIENCY

In papers [11] and [12] Fisher sketched his influential program of establishing the asymptotic efficiency of maximum likelihood estimators. (See [1, 44, 45] for background.) Let $\hat{\theta}_n$ be the maximum likelihood estimate for a scalar parameter $\theta$ found from a sample of size $n$. Fisher implicitly assumed regularity conditions that implied the existence of maximum likelihood estimates and much more. In the general discussion of this section we will not mention explicitly the required regularity conditions.

Fisher's idea was to prove that:

1. The scaled difference $(\hat{\theta}_n - \theta)n^{1/2}$ is asymptotically normal with parameters $(0, 1/I(\theta))$, where $I(\theta)$ is Fisher's information. [In this paper the normal distribution $N(\mu, \sigma^2)$ is parameterized by its expectation $\mu$ and its variance $\sigma^2$.]

2. If another estimator $T_n$ is such that $(T_n - \theta)n^{1/2}$ is asymptotically normal with parameters $(0, v(\theta))$, then $v(\theta) \geq 1/I(\theta)$.

Fisher proposed several informal arguments for these two statements. The first statement was established rigorously by Cramér [8]. Later Cramér's regularity conditions were relaxed, and analogous statements were established for methods of estimation different from maximum likelihood (such as Bayes estimators or Weiss and Wolfowitz's [51] maximum probability estimators). The second statement is wrong if understood literally, as shown by Hotelling in his letter to Fisher (see the epigraph; available on-line [17] and quoted by Stigler in [45]).

The bluntest interpretation of Hotelling's objection is that, for each parameter value $\theta$, the estimator that is identically equal to $\theta$,

$$T_n := \theta, \tag{1}$$

is such that $(T_n - \theta)n^{1/2}$ is asymptotically normal with parameters $(0, 0)$. Since $0 < 1/I(\theta)$, the parameter point $\theta$ will be a "point of superefficiency." This notion of superefficiency was perhaps not particularly interesting to Fisher and Hotelling, since the estimator (1) is not consistent at parameter points different from $\theta$. Hodges's implementation of Hotelling's idea (probably discovered completely independently) is to set

$$T_n := \begin{cases} \hat{\theta}_n, & \text{if } |\hat{\theta}_n - \theta| \geq n^{-1/4}, \\ \theta, & \text{if } |\hat{\theta}_n - \theta| < n^{-1/4} \end{cases} \tag{2}$$

[Le Cam [26], Section 1, with a credit to Hodges (1951); Le Cam says that Hodges produced a series of examples and gives an example slightly different from (2)]. The advantage of Hodges's estimator is that it is consistent and, moreover, its asymptotic expected squared error is never worse than that of the maximum likelihood estimator (at least in the case of the Gaussian model with variance 1 considered by Le Cam). Hodges's estimator may be said to be superefficient at $\theta$ in the narrow sense: asymptotically, it beats the maximum likelihood estimator at $\theta$ *and* is not worse than the maximum likelihood estimator at the other parameter points. The estimator (1) is then superefficient in the wide sense.

We will refer to the three approaches to dealing with superefficiency, the lines of defense mentioned in Section 1, as the first approach (exclude the evil by changing the qualifying rules), the second approach (show that the evil, i.e., the set of points of superefficiency, is not great), and the third approach (declare a parameter point a point of inefficiency if some of its neighbors are points of inefficiency). This appears to be the chronological order of their appearance. Some work in broadly the same direction, such as that on the Bahadur [2] and Rao ([38], Definitions 2.3–2.6) efficiency of estimators, is of a very different character and cannot be easily assigned to one of the three approaches.

## 2.1 Exclusion of the Evil

It appears that the first approach was initiated by Fisher himself in 1930, who, in response to Hotelling's doubts, gave his "third proof" of the efficiency of the maximum likelihood estimator (in the terminology of Stigler [45], who points out that of Fisher's three proofs this is the only real proof). Fisher considered only a finite observation space and restricted competition by considering only consistent estimators that are smooth functions (independent of the sample size $n$) of the observed relative frequencies $x_i$. A modification of the "third proof" was published in [13] (pages 45–46), which considered the consistent estimators defined by an equation of the form $\sum_i k_i(\theta)x_i = 0$, the summation being over the observation space.

A simple proof of Fisher's bound for finite observation spaces and consistent estimators that are smooth functions of the observed relative frequencies was given by Rao in 1955 [37] (and reproduced in [40]). This proof was extended by Kallianpur and Rao [21] to the observation space $\mathbb{R}$; they considered estimators that are Fréchet differentiable functions of the empirical distribution function. Fréchet differentiability was weakened to Gâteaux differentiability by Kallianpur [20].

Another restricted class of estimators (although more general than Fisher's) was considered by Neyman [33]. Neyman's overview of known properties of maximum likelihood estimators reflects beliefs prevailing at the time. Wolfowitz starts his review of [33] in *Mathematical Reviews* as follows:

> It is well known that maximum likelihood (ML) estimates have, under general conditions, the following properties: (a) consistency, (b) asymptotic normality, (c) minimal variance of the limiting distribution.

The corresponding statement in Neyman's paper is more hedged; Neyman refers to the earlier work by Hotelling [18] and Doob [9], neither of whom, however, discussed (c).

An important byproduct of the work on the first approach was the Cramér–Rao inequality for unbiased estimators ([14, 36], [8], Section 32.3). It appears that this result *per se* is not directly connected with Fisher's program (as emphasized by Weiss and Wolfowitz, [51], pages 10–11). As a consequence, results about superefficiency that are based on the Cramér–Rao inequality (such as Theorem 1 in [50]) impose regularity conditions on the allowed estimators that are difficult to interpret.

The first approach has often been criticized. For example, Wolfowitz ([52], page 249) writes:

> ...to argue that the maximum likelihood (m.l.) estimator is best by ruling out some of its competitors, is a dangerous if tempting procedure. It can easily result in begging the entire question. After all, to give an example from social life, anyone can become the chess champion of his town if the better players are arbitrarily declared ineligible to compete. Yet what we are seeking to establish is that the m.l. estimator is asymptotically the champion!

In particular, he objects against the assumption of asymptotical normality of the estimators admitted to the competition. This is echoed by Weiss and Wolfowitz [51]:

> The problem is, however, to exclude *only* artificial competitors. If we exclude sensible and practical competitors then any claims about the optimality of the m.l. or any other estimator are hollow indeed, and the theorems proved do not describe the physical reality and are not of practical value or aesthetic interest.

In Wolfowitz's [52] terminology, any regularity conditions imposed on the estimators should be "statistically operational." He believed that the weak uniform convergence of $(T_n - \theta)n^{1/2}$ to a random variable (not necessarily normal) depending on $\theta$ is such a statistically operational condition. The requirement of weak uniform convergence was also proposed by Rao [39] in 1963 (the same year that the results of [52] were presented at the Seventh All–Soviet Union Conference on Probability and Mathematical Statistics). Lehmann [29] suggests the alternative condition that the variance $v(\theta)$ of the limiting distribution of $(T_n - \theta)n^{1/2}$ should be a continuous function of $\theta$. Lehmann notices that his condition is weaker than the condition of weak uniform convergence (under mild regularity conditions on the statistical model; cf. [52], Lemma 2) but also eliminates superefficiency: this follows immediately from Le Cam's result, since superefficiency at one point leads to superefficiency in a neighborhood of that point when $v$ is continuous.

Pfanzagl [34] develops further Wolfowitz's objection against Fisher's program:

With the same justification with which Wolfowitz questioned the asymptotic normality assumption for the sequence of estimates one could question his assumption of weak uniform convergence: Why should a statistician confine himself to estimates for which the sequence of distributions of $n^{1/2}(T_n - \theta)$ converges at all?

In his Theorem 1 Pfanzagl proves the absence of points of superefficiency for median unbiased estimators ([34], Theorem 1); this result is extended by Michel [32] to what he calls strongly asymptotically median unbiased estimators.

## 2.2 Deprecation of the Evil

The second approach was started by Le Cam's result that sets of points of superefficiency have Lebesgue measure zero. The earliest version of this result was given without proof in the abstract [25] of his dissertation [26]. The dissertation itself [26] states and contains a proof of a stronger version of the result (cf. [25], Theorem 2, and [26], Theorem 9). However, the proof given in [26] is wrong, as noticed by Wolfowitz [52]; it does not even prove the weaker version of [25]. Corrected proofs were given by Le Cam himself [27] and Bahadur [3]. Paper [48] is devoted to the history and several proofs of Le Cam's result.

The main difference between the versions of Le Cam's result given in [25] and [26] is that [26] does not assume that the estimator $T_n$ is asymptotically normal, whereas [25] makes this assumption. Le Cam [27] and Bahadur [3] revert to asymptotically normal estimators. Pfanzagl ([34], Theorem 2) removes the condition of asymptotic normality proving a result similar to the one claimed in Le Cam [26]. Both Le Cam [26] and Pfanzagl [34] assume that some function of $T_n$ ($n(T_n - \theta)^2$ in [26] and $(T_n - \theta)n^{1/2}$ in [34]) converges weakly to some probability measure. Therefore, all these papers involve elements of the first approach.

Whereas Le Cam [25, 26] considers superefficiency in the narrow sense, the results given in [3, 27, 34] concern superefficiency in an intermediate sense: the assumptions made about the estimator $T_n$ imply its consistency (and more), but it is not required that the asymptotic variance of $(T_n - \theta)n^{1/2}$ should never exceed $1/I(\theta)$.

In his paper [26] Le Cam claims that sets of points of superefficiency can be uncountable. There is no formal contradiction between Le Cam's claim and this paper's result: in his example (Example 4 in [26]) Le Cam uses a different, somewhat arbitrary, notion of superefficiency. This example will be further discussed in Section 3.3.

The standard textbook [7], page 305, asserts that sets of points of superefficiency are countable. However, this is simply a slip of the pen, since this statement is attributed to Le Cam [26], who never makes it.

## 2.3 Collective Responsibility

In the third approach, when evaluating the performance of an estimator $T_n$ at a parameter point $\theta$, one takes into account the performance of $T_n$ at parameter points different from $\theta$. As discussed at the beginning of this section, there is a whiff of this already in the standard notion of superefficiency ([26], Definition 4), as used in the Berkeley group in the early 1950s: $\theta$ does not qualify as a point of superefficiency of (1) because $T_n$ is so inefficient, not even consistent, at all other points.

In the last section of his paper [26] Le Cam proves several results that belong to the third approach. His Theorem 14 says that the performance of a superefficient estimator in a shrinking neighborhood of a point of superefficiency is poor. His Theorem 13 states this result in terms of a formal measure of performance of an estimator taking into account the performance at the neighboring points.

Another early paper explicitly using the third approach is Chernoff's [5]. Theorem 1 of that paper, in Chernoff's words, "states that for an arbitrary estimate *the reciprocal of the information is 'essentially' asymptotically a lower bound for the asymptotic variance.*" The word "essentially" refers to taking the supremum of the asymptotic variances (suitably modified) over a shrinking neighborhood of the given parameter point.

The culmination of this line of work was Hájek's [16] local asymptotic minimax theorem. (See Le Cam [28], pages 24–25, for a discussion of connections of this theorem with other results.) Hájek's result has been generalized in various directions, and at this time the third approach is perhaps the dominant one.

## 2.4 Informal Comparison

The difference between the first and third approaches is not always clear-cut. If an estimator performs well at a parameter point $\theta$ but much worse at $\theta$'s neighbors, we can react to this in two ways: either eliminate the estimator from competition (first approach) or punish the estimator by declaring its performance at $\theta$ to be its worst performance at $\theta$ and its neighbors (third approach). The former option is implemented as, for example, the requirement of weak uniform convergence in [39, 52] (discussed in Section 2.1), the requirement of continuous convergence in [42] and requirement (3.7) in [51]. Apart from this borderline situation, the objections against the first approach quoted in Section 2.1 appear to be valid.

The second and third approaches are convincing in different circumstances. The third approach is convincing when our *a priori* expectations for various values of $\theta$ are diffuse. In the Bayesian case, where these expectations are expressed via a full-blown prior distribution, this distribution should not assign a positive weight to any specific value of $\theta$. If the expectations are not diffuse (e.g., when the value $\theta = 0$ corresponds to no difference between two treatments, the statistician might want to assign to it a positive probability), or too uncertain for us to judge how diffuse they are, the third approach becomes less convincing.

### 2.5 Contribution of this Paper

Our main result, Theorem 1 in the next section, answers a natural question: can a given parameter point $\theta$ be a point of superefficiency? Hotelling's and Hodges's examples, (1) and (2), work for any $\theta$, but if $\theta$ is noncomputable, the resulting estimators are also noncomputable, and their existence is of no use. Our theorem says that no computable estimator can be superefficient at a noncomputable point. In this way the notion of computability further limits the damage inflicted by Hotelling's objection: yes, superefficiency (in its most extreme form, $T_n = \theta$ for all $n$) is possible at computable points $\theta$, but there can be no superefficiency at the other $\theta$.

## 3. COMPUTABILITY OF POINTS OF SUPEREFFICIENCY FOR COMPUTABLE ESTIMATORS

Let $\Omega_1, \Omega_2, \ldots$ be a sequence of measurable spaces, and for each $n \in \{1, 2, \ldots\}$, let $\{P_{n,\theta} \mid \theta \in \Theta\}$ be a statistical model on $\Omega_n$. We will assume that $\Theta$ is an open interval of the real line ($\Theta = \mathbb{R}$ is allowed). Each $P_{n,\theta}$ is a probability measure on $\Omega_n$, and $\theta \in \Theta$ is the parameter to be estimated. Little will be lost if the reader assumes that $\Omega_n = \Omega^n$ and $P_{n,\theta} = (P_\theta)^n$ for all $n$, which corresponds to independent observations chosen from an observation space $\Omega$ according to $P_\theta$. An *estimator* $T = \{T_n\}_{n=1}^\infty$ for $\{P_{n,\theta}\}$ is a sequence of measurable functions $T_n : \Omega_n \to \Theta$.

We will need a condition of regularity, which will be stated in terms of a natural measure of closeness between probability measures. Formally, the *affinity* between probability measures $P$ and $Q$ on the same measurable space $\Omega$ is defined by

$$(3) \qquad \pi(P, Q) := \inf_E \max(P(E), Q(\Omega \setminus E)),$$

$E$ ranging over the measurable sets in $\Omega$. Notice that:

- it is always true that $\pi(P, Q) \in [0, 1]$ and $\pi(P, P) \in [1/2, 1]$;
- probability measures $P$ and $Q$ are mutually singular if and only if $\pi(P, Q) = 0$;
- sequences of probability measures $P_n$ and $Q_n$ on measurable spaces $\Omega_n$ are asymptotically entirely separated if and only if $\liminf_{n \to \infty} \pi(P_n, Q_n) = 0$.

Another ingredient of our regularity condition will be a continuous function $I : \Theta \to (0, \infty)$ (typically, Fisher's information).

ASSUMPTION 1. *For any $\varepsilon > 0$ and any $\theta \in \Theta$, there exist a positive integer $N$ and a neighborhood $O \subseteq \Theta$ of $\theta$ such that, for all $n \geq N$ and $\theta_1, \theta_2 \in O$,*

$$(4) \qquad \pi(P_{n,\theta_1}, P_{n,\theta_2}) \geq \Phi\left(-\frac{|\theta_2 - \theta_1|\sqrt{nI(\theta)}}{2}\right) - \varepsilon,$$

*where $\Phi$ is the $N(0, 1)$ distribution function.*

Assumption 1 is a weak form of the uniform condition of local asymptotic normality. It will be discussed in the next section. In Section 4.2 we will see that it is satisfied for statistical models satisfying standard regularity conditions (cf. the reference to [19] there). As a simple sanity test, in Sections 4.1 and 4.3 we check it for the Gaussian model with known variance.

THEOREM 1. *Suppose $\{T_n\}$ is a computable estimator for $\{P_{n,\theta}\}$ satisfying Assumption 1. For any $c > 0$ and any noncomputable $\theta \in \Theta$,*

$$(5) \qquad \limsup_{n \to \infty} P_{n,\theta}(|T_n - \theta| > cn^{-1/2}) \geq \Phi(-c\sqrt{I(\theta)}).$$

As we said earlier, computability is discussed in Appendix A. The reader who is only interested in the countability of points of superefficiency (Corollary 2 below) can ignore all statements about computability; the proof of Theorem 1 will still show that there are only countably many points of super-efficiency under Assumption 1. A streamlined independent proof of the countability of points of superefficiency is given in Appendix B.

If $(\hat{\theta}_n - \theta)n^{1/2}$ is asymptotically $N(0, 1/I(\theta))$ under $P_{n,\theta}$ for some estimator $\hat{\theta}_n$, such as the maximum likelihood estimator $\hat{\theta}_n$, we will have an "almost opposite" inequality to (5),

$$
\begin{aligned}
(6) \quad & \limsup_{n\to\infty} P_{n,\theta}(|\hat{\theta}_n - \theta| > cn^{-1/2}) \\
& \leq 2\Phi(-c\sqrt{I(\theta)}).
\end{aligned}
$$

The use of probabilities $P_{n,\theta}(|T_n - \theta| > cn^{-1/2})$ for measuring the concentration of estimators is very standard in the literature discussed in Section 2: cf., for example, Le Cam's discussion of concentration in [26] (starting from page 288), Wolfowitz's Theorem [52] (pages 258–259), Schmetterer's [42] Theorem 2.2, Pfanzagl's [34] Theorems 1 and 2.

There is a gap between the right-hand sides of (5) and (6). To eliminate it, we can consider only large values of $c$. Let us define the *asymptotic efficiency* of an estimator $T = \{T_n\}$ at a parameter point $\theta \in \Theta$ as

$$
\begin{aligned}
(7) \quad & \mathrm{ae}_\theta(T) \\
& := \liminf_{c\to\infty} \liminf_{n\to\infty} \frac{-\ln P_{n,\theta}(|T_n - \theta| > cn^{-1/2})}{c^2 I(\theta)/2}
\end{aligned}
$$

(with convention $-\ln 0 := \infty$). Since $-\ln\Phi(-x) \sim x^2/2$ as $x \to \infty$ (see, e.g., [10], Lemma VII.2), (6) implies that $\mathrm{ae}_\theta(\hat{\theta}) \geq 1$. In this sense the maximum likelihood estimators are efficient under the usual regularity conditions. We can say that $T$ is super-efficient at $\theta$ if $\mathrm{ae}_\theta(T) > 1$. Inequality (5) implies $\mathrm{ae}_\theta(T) \leq 1$.

In the classical case of $(T_n - \theta)n^{1/2}$ asymptotically normal with parameters $(0, v(\theta))$, $\mathrm{ae}_\theta(T) = 1/(I(\theta)v(\theta))$, and so, under the usual regularity conditions, $\mathrm{ae}_\theta(T)$ is the ratio of the asymptotic variance of the rescaled maximum likelihood estimator to the asymptotic variance of rescaled $T_n$. Therefore, in this case $\mathrm{ae}_\theta(T)$ is the classical asymptotic efficiency of $T$ at $\theta$, as defined by Fisher [11] (page 316) and Cramér [8] (Section 32.5).

Before proving Theorem 1, we state three simple corollaries of it, all assuming Assumption 1.

COROLLARY 1. *If the parameter point $\theta$ is non-computable and a computable estimator $\{T_n\}$ is such that $(T_n - \theta)n^{1/2}$ weakly converges to $N(0, v(\theta))$, then $v(\theta) \geq 1/I(\theta)$.*

PROOF. It suffices to notice that $1/(I(\theta)v(\theta)) = \mathrm{ae}_\theta(T) \leq 1$.  □

COROLLARY 2. *If $c > 0$ and $\{T_n\}$ is an estimator for $\{P_{n,\theta}\}$, the inequality*

$$
(8) \quad \limsup_{n\to\infty} P_{n,\theta}(|T_n - \theta| > cn^{-1/2}) < \Phi(-c\sqrt{I(\theta)})
$$

*holds for at most countably many $\theta \in \Theta$. In particular, $\mathrm{ae}_\theta(T) > 1$ for at most countably many $\theta$. In particular, if $(T_n - \theta)n^{1/2}$ weakly converges to $N(0, v(\theta))$ for all $\theta \in \Theta$, then $v(\theta) < 1/I(\theta)$ holds for at most countably many $\theta$.*

The last part of Corollary 2 was also proved (under different regularity conditions) in [49], the Appendix.

PROOF OF COROLLARY 2. Every estimator is computable with respect to some oracle (see Section A.4 of Appendix A for information about oracles). Fix such an oracle for $\{T_n\}$. Theorem 1 will continue to hold if computability is replaced by computability with respect to this oracle. Finally, there are only countably many parameter points that are computable with respect to this oracle.  □

As mentioned earlier, a proof of Corollary 2 not using the notions of computability and oracle can be extracted from the proof of Theorem 1. See Appendix B for details.

We can define the *asymptotic estimability* $\mathrm{ae}_\theta$ of a parameter point $\theta$ as

$$
\mathrm{ae}_\theta := \sup_T \mathrm{ae}_\theta(T),
$$

with $T$ ranging over the computable estimators. The following corollary is a formalization of a new all-or-nothing phenomenon arising in our algorithmic framework.

COROLLARY 3. *Suppose there is a computable estimator $\hat{\theta}_n$ satisfying (6). Then, for each $\theta \in \Theta$, either $\mathrm{ae}_\theta = 1$ or $\mathrm{ae}_\theta = \infty$.*

PROOF. By (6) we have $\mathrm{ae}_\theta \geq 1$. In combination with (5), this gives $\mathrm{ae}_\theta = 1$ for noncomputable $\theta$. If $\theta$ is computable, setting (1) gives $\mathrm{ae}_\theta = \infty$.  □

The formalization given by Corollary 3 is imperfect, because we have much more than $\mathrm{ae}_\theta = \infty$ for computable $\theta$: there is a computable estimator that estimates $\theta$ with zero error, $T_n = \theta$.

### 3.1 The Role of Computability

In Theorem 1 and its corollaries the notion of computability is applied to two kinds of objects: estimators and parameter points. There is, however, an important philosophical difference between computable estimators and computable parameter points.

The purpose of estimators is to be used for computing estimates, and so their computability is essential. Accordingly, in our discussion we restrict ourselves to computable estimators.

A parameter point is not meant to be computed by anybody. Depending on which school of statistics we listen to, it is either a constant chosen by Nature or a mathematical fiction. In any case, there is no reason for us to require or expect that it should be computable. In fact, noncomputable parameter points are often more important than computable ones: for example, if the parameter point is chosen from a smooth prior on the real line, it will be noncomputable with probability one.

Theorem 1 implies that, for the standard regular statistical models, the maximum likelihood estimator is efficient (cannot be beaten by any other computable estimator) at noncomputable parameter points. This statement remains important despite the complementary statement that the maximum likelihood estimator is greatly outperformed by Hotelling's and Hodges's estimators at computable parameter points.

### 3.2 Proof of Theorem 1

The proof will use the following implication of Assumption 1.

LEMMA 1. *Let $c > 0$, $a \in (0, 1)$, $\theta \in \Theta$, and $I \geq I(\theta)$. There exist $\varepsilon > 0$, positive integer $N$, and a neighborhood $O \subseteq \Theta$ of $\theta$ such that, for any $n \geq N$, $\theta_1, \theta_2 \in O$, and $A_1, A_2 \subseteq \Omega_n$ satisfying*

$$(9) \qquad \max(P_{n,\theta_1}(A_1), P_{n,\theta_2}(A_2)) \leq a\Phi(-c\sqrt{I})$$

*and*

$$(10) \qquad |\theta_2 - \theta_1| \leq 2(1+\varepsilon)^3 cn^{-1/2},$$

*it is true that $A_1 \cup A_2 \neq \Omega_n$.*

PROOF. Let $\varepsilon > 0$ be so small that

$$\Phi(-(1+\varepsilon)^3 c\sqrt{I}) - \varepsilon > a\Phi(-c\sqrt{I}).$$

Take any $N$ and $O$ satisfying the condition in Assumption 1. Using $I(\theta) \leq I$ and (10), we now obtain,

for $n \geq N$,

$$\pi(P_{n,\theta_1}, P_{n,\theta_2}) \geq \Phi\left(-\frac{|\theta_2 - \theta_1|\sqrt{nI(\theta)}}{2}\right) - \varepsilon$$

$$\geq \Phi(-(1+\varepsilon)^3 c\sqrt{I}) - \varepsilon$$

$$> a\Phi(-c\sqrt{I}).$$

Were it true that $A_1 \cup A_2 = \Omega_n$, (9) would imply the opposite inequality

$$\pi(P_{n,\theta_1}, P_{n,\theta_2}) \leq \max(P_{n,\theta_1}(A_1), P_{n,\theta_2}(\Omega \setminus A_1))$$

$$\leq \max(P_{n,\theta_1}(A_1), P_{n,\theta_2}(A_2))$$

$$\leq a\Phi(-c\sqrt{I}). \qquad \square$$

PROOF OF THEOREM 1. Fix some $c > 0$ and $\theta \in \Theta$ such that (5) fails. We will exhibit an algorithm for computing $\theta$, which will prove the theorem. There exist $a \in (0, 1)$, $I > I(\theta)$, and $N$ such that, for all $n \geq N$,

$$(11) \qquad P_{n,\theta}(|T_n - \theta| > cn^{-1/2}) \leq a\Phi(-c\sqrt{I}).$$

Let $\varepsilon > 0$, $N$, and $O \ni \theta$ satisfy the condition in Lemma 1 [so that $N$ is assumed to be large enough both for (11) to hold and for the condition in Lemma 1 to be satisfied].

Choose an open interval $(L, R) \subseteq O$ with rational end-points such that $\theta \in (L, R)$ and $I(q) < I$ for all $q \in (L, R)$. In what follows we will also impose some other conditions on the interval $(L, R)$ (it has to be sufficiently short).

Let $(\theta_1, \theta_2) \subseteq (L, R)$ be an open interval containing $\theta$. We will construct, in a computable manner, an open interval $(\theta_1', \theta_2') \subseteq (L, R)$ whose length is at most $(1+\varepsilon)^{-1}|\theta_2 - \theta_1|$ and which still contains $\theta$. Repeating this operation, we can compute $\theta$ to any accuracy starting from $(L, R)$: the length of the interval known to contain $\theta$ will tend to zero exponentially fast.

First notice that we can compute a positive integer $n$ such that

$$(12) \qquad \begin{aligned} 2(1+\varepsilon)^2 cn^{-1/2} &\leq |\theta_2 - \theta_1| \\ &\leq 2(1+\varepsilon)^3 cn^{-1/2}, \end{aligned}$$

and that we can assume that the resulting $n$ satisfies $n \geq N$. Indeed, the double inequality (12) is equivalent to

$$(13) \qquad \frac{4(1+\varepsilon)^4 c^2}{(\theta_2 - \theta_1)^2} \leq n \leq \frac{4(1+\varepsilon)^6 c^2}{(\theta_2 - \theta_1)^2},$$

so making $(L, R) \supseteq (\theta_1, \theta_2)$ sufficiently short will ensure the existence of $n$ satisfying (12) and the inequality $n \geq N$ for all such $n$.

Let us say that a point $q \in \Theta$ is *suitable* if

(14) $\qquad P_{n,q}(|T_n - q| > cn^{-1/2}) \leq a\Phi(-c\sqrt{I})$.

Let $\mathfrak{S} = \mathfrak{S}_{n,\theta_1,\theta_2}$ be the set of all suitable points in $(\theta_1, \theta_2)$. We know that $\theta \in \mathfrak{S}$; see (11).

Let us show that $|q_2 - q_1| \leq (1+\varepsilon)^{-2}|\theta_2 - \theta_1|$ for all $q_1, q_2 \in \mathfrak{S}$. This follows from Lemma 1: setting $A_i := \{\omega \mid |T_n(\omega) - q_i| > cn^{-1/2}\}$, $i = 1, 2$, and using (12) and (14), we can see that there exists $\omega \in \Omega_n \setminus (A_1 \cup A_2)$. Since

$$|T_n(\omega) - q_1| \leq cn^{-1/2}$$

and

$$|T_n(\omega) - q_2| \leq cn^{-1/2},$$

the triangle inequality and (12) imply

$$|q_2 - q_1| \leq 2cn^{-1/2} \leq (1+\varepsilon)^{-2}|\theta_2 - \theta_1|.$$

The estimator $\{T_n\}$ is computable, and so we can compute an open interval $(\theta_1', \theta_2') \supseteq \mathfrak{S}$ of length $|\theta_2' - \theta_1'| \leq (1+\varepsilon)^{-1}|\theta_2 - \theta_1|$. (See Section A.2 of Appendix A for details.) This completes the proof of the theorem. □

### 3.3 Le Cam's Example

In one of the examples in [26] (Example 4), Le Cam constructs an uncountable set $S$ and an estimator $\{T_n\}$ for the Gaussian model $\{P_{n,\theta}\} = \{(N(\theta, 1))^n\}$ with unknown mean $\theta$ and known variance 1 such that: (a) for all $\theta \in S$ and for all $n$ of the form $n = 7^{2j}$, $j = 2, 3, \ldots$,

$$P_{n,\theta}(|T_n - \theta| > n^{-1/2})$$
$$\leq P_{n,\theta}(|\hat{\theta}_n - \theta| > n^{-1/2}) - 0.18,$$

and (b) for all $\theta \notin S$ and for all sufficiently large $n$ of the form $n = 7^{2j}$,

$$P_{n,\theta}(|T_n - \theta| > n^{-1/2}) = P_{n,\theta}(|\hat{\theta}_n - \theta| > n^{-1/2}),$$

$\hat{\theta}_n$ being the maximum likelihood estimator. This does not contradict our Corollary 2 for two reasons: first, our definition of $\mathrm{ae}_\theta(T)$ involves the probabilities $P_{n,\theta}(|T_n - \theta| > cn^{-1/2})$ for large $c$, whereas Le Cam arbitrarily fixes $c := 1$, and second, the restriction to $n = 7^{2j}$ invalidates our proof, which crucially depends on the set of allowed $n$ being sufficiently dense [see (13)].

At the end of his construction of $\{T_n\}$, Le Cam says that a different, "necessarily more complicated," construction will give another estimator and another uncountable set that enjoy similar properties to $\{T_n\}$ and $S$ without the restriction to $n = 7^{2j}$. This would eliminate the second reason, but an independent verification of Le Cam's claim is desirable: there is no proof in [26], and even the proofs in that paper are "quelquefois incorrectes" ([27], page 17).

## 4. REGULARITY CONDITIONS

The present section and Appendix A are devoted to the regularity conditions imposed on the sequence of statistical models $\{P_{n,\theta}\}$ and the estimator $\{T_n\}$, respectively. The status of these two sets of regularity conditions is very different: whereas the conditions imposed on the estimator should be minimal (cf. Section 2.1), we can be much more flexible in choosing the conditions imposed on the sequence of statistical models: even if these conditions are relatively strong, they are still likely to be satisfied by many important models (cf. the discussion in [52], Section 3).

In this section we will see that Assumption 1, essentially saying that

$$\liminf_{\substack{\theta_1, \theta_2 \to \theta \\ n \to \infty}} \left(\pi(P_{n,\theta_1}, P_{n,\theta_2}) - \Phi\left(-\frac{|\theta_2 - \theta_1|\sqrt{nI(\theta)}}{2}\right)\right)$$
$$\geq 0$$

for all $\theta \in \Theta$, follows from easier to check or more standard assumptions.

### 4.1 In Terms of the Likelihood Ratio

We make the standard assumption that for each $n$ all $P_{n,\theta}$ are absolutely continuous with respect to a $\sigma$-finite measure $\mu_n$. Let $f_{n,\theta}$ be a density of $P_{n,\theta}$ with respect to $\mu_n$.

ASSUMPTION 2. *For any $\varepsilon > 0$ and any $\theta \in \Theta$, there exist a positive integer $N$ and a neighborhood $O \subseteq \Theta$ of $\theta$ such that, for all $n \geq N$ and for all distinct $\theta_1, \theta_2 \in O$,*

(15)
$$P_{n,\theta_1}\left(\frac{f_{n,\theta_2}}{f_{n,\theta_1}} > 1\right)$$
$$\geq \Phi\left(-\frac{|\theta_2 - \theta_1|\sqrt{nI(\theta)}}{2}\right) - \varepsilon.$$

LEMMA 2. *Assumption 2 implies Assumption 1.*

PROOF. By symmetry, we can complement (15) by

$$P_{n,\theta_2}\left(\frac{f_{n,\theta_1}}{f_{n,\theta_2}} > 1\right) \geq \Phi\left(-\frac{|\theta_1 - \theta_2|\sqrt{nI(\theta)}}{2}\right) - \varepsilon.$$

Therefore, it suffices to prove

$$(16) \quad \pi(P_{n,\theta_1}, P_{n,\theta_2}) \geq \min\left(P_{n,\theta_1}\left(\frac{f_{n,\theta_2}}{f_{n,\theta_1}} > 1\right),\right.$$
$$\left. P_{n,\theta_2}\left(\frac{f_{n,\theta_1}}{f_{n,\theta_2}} > 1\right)\right).$$

Suppose this inequality is false, and so we can find $t$ such that

$$\pi(P_{n,\theta_1}, P_{n,\theta_2})$$
$$< t$$
$$< \min\left(P_{n,\theta_1}\left(\frac{f_{n,\theta_2}}{f_{n,\theta_1}} > 1\right), P_{n,\theta_2}\left(\frac{f_{n,\theta_1}}{f_{n,\theta_2}} > 1\right)\right).$$

In this case, there exists an event $E$ such that

$$P_{n,\theta_1}(E) < t, \quad P_{n,\theta_1}\left(\frac{f_{n,\theta_2}}{f_{n,\theta_1}} > 1\right) > t,$$
$$P_{n,\theta_2}(E) > 1 - t, \quad P_{n,\theta_2}\left(\frac{f_{n,\theta_2}}{f_{n,\theta_1}} \geq 1\right) < 1 - t.$$

This contradicts the Neyman–Pearson lemma. $\square$

An easy calculation shows that Assumption 2 is satisfied for sampling from the Gaussian family $N(\theta, \sigma^2)$ with known variance $\sigma^2 > 0$ and with $I(\theta) := \sigma^{-2}$ [in other words, with $I(\theta)$ Fisher's information]. In fact, for this statistical model we have

$$P_{n,\theta_1}\left(\frac{f_{n,\theta_2}}{f_{n,\theta_1}} > 1\right) = \Phi\left(-\frac{|\theta_2 - \theta_1|\sqrt{n}}{2\sigma}\right),$$

with an equality and without the need to subtract $\varepsilon$. Therefore, Assumption 2 and, *a fortiori*, Assumption 1 are not vacuous.

## 4.2 Local Asymptotic Normality

Another assumption that implies Assumption 1 is the following uniform version of the condition of local asymptotic normality.

ASSUMPTION 3. For any $\theta \in \Theta$, any $\lambda \geq 0$, any sequence $\theta_i \to \theta$ of elements of $\Theta$, any sequence $n_i \to \infty$ of positive integers, and any sequence $\lambda_i \to \lambda$ of positive real numbers such that $\theta_i + \lambda_i/\sqrt{n_i I(\theta_i)} \in \Theta$ for all $i = 1, 2, \ldots$, there exist sequences $\Delta_i$ and $\psi_i$ of random variables such that, for all $i$,

$$(17) \quad \ln \frac{f_{n_i, \theta_i + \lambda_i/\sqrt{n_i I(\theta_i)}}}{f_{n_i, \theta_i}} = \lambda \Delta_i - \lambda^2/2 + \psi_i$$

and:

- the distribution of $\Delta_i$ with respect to $P_{n_i, \theta_i}$ weakly converges to $N(0, 1)$;
- $\psi_i$ converges to 0 in $P_{n_i, \theta_i}$-probability.

The derivation of a slightly stronger version of Assumption 3 under standard regularity conditions can be found in [19] (Definition II.2.2, Theorem II.1.2 and Remark II.1.4).

LEMMA 3. *Assumption 3 implies Assumption 1.*

PROOF. Suppose that Assumption 3 holds whereas Assumption 1 does not hold. The latter implies that there exist $\varepsilon > 0$ and $\theta \in \Theta$ such that for each positive integer $N$ and each neighborhood $O$ of $\theta$ there exist $n \geq N$ and $\theta_1, \theta_2 \in O$ for which (4) is violated. Fix such $\varepsilon$ and $\theta$. There exist sequences $n_i \to \infty$ of positive integers and $\theta_i \to \theta$, $\tilde{\theta}_i \to \theta$ of elements of $\Theta$ such that, for all $i$, $\theta_i < \tilde{\theta}_i$ and

$$(18) \quad \pi(P_{n_i, \theta_i}, P_{n_i, \tilde{\theta}_i}) < \Phi\left(-\frac{(\tilde{\theta}_i - \theta_i)\sqrt{nI(\theta)}}{2}\right) - \varepsilon.$$

It is clear that $\theta$ in (18) can be replaced by $\theta_i$ (slightly decrease $\varepsilon$ and disregard the initial $i$'s if necessary). Setting $\lambda_i := (\tilde{\theta}_i - \theta_i)\sqrt{n_i I(\theta_i)}$, we can rewrite (18) as

$$(19) \quad \pi(P_{n_i, \theta_i}, P_{n_i, \tilde{\theta}_i}) < \Phi\left(-\frac{\lambda_i}{2}\right) - \varepsilon.$$

The last inequality shows that the sequence $\lambda_i$ is bounded. Therefore, we can assume, without loss of generality, that $\lambda_i \to \lambda$ for some $\lambda \geq 0$ (consider a subsequence of $i$ if necessary). Fix sequences $\Delta_i$ and $\psi_i$ satisfying the conditions in Assumption 3.

Notice that $\lambda > 0$: indeed, if $\lambda$ were zero, (17) would converge to zero in probability, which would contradict (19). Therefore,

$$(20) \quad P_{n_i, \theta_i}\left(\frac{f_{n_i, \tilde{\theta}_i}}{f_{n_i, \theta_i}} > 1\right) = P_{n_i, \theta_i}\left(\Delta_i > \frac{\lambda}{2} - \frac{\psi_i}{\lambda}\right)$$
$$\to \Phi\left(-\frac{\lambda}{2}\right).$$

In a similar way we can obtain

$$(21) \quad P_{n_i, \tilde{\theta}_i}\left(\frac{f_{n_i, \theta_i}}{f_{n_i, \tilde{\theta}_i}} > 1\right) \to \Phi\left(-\frac{\lambda}{2}\right).$$

Inequalities (20) and (21) contradict (19) and (16). $\square$

### 4.3 In Terms of the Variation Distance

The *variation distance* $\|P - Q\|$ between two probability measures on the same measurable space $\Omega$ is defined to be

$$\|P - Q\| := \sup_E |P(E) - Q(E)|,$$

$E$ ranging over the measurable sets in $\Omega$. A slightly stronger form of Assumption 1 can be stated in terms of variation distance rather than affinity.

ASSUMPTION 4. For any $\varepsilon > 0$ and any $\theta \in \Theta$, there exist a positive integer $N$ and a neighborhood $O \subseteq \Theta$ of $\theta$ such that, for all $n \geq N$ and $\theta_1, \theta_2 \in O$,

$$\|P_{n,\theta_1} - P_{n,\theta_2}\| \leq 1 - 2\Phi\left(-\frac{|\theta_2 - \theta_1|\sqrt{nI(\theta)}}{2}\right) + \varepsilon.$$

Theorem 1 remains true if Assumption 1 is replaced by Assumption 4. This follows from the following lemma.

LEMMA 4. *It is always true that*

$$\pi(P, Q) \geq \frac{1 - \|P - Q\|}{2}.$$

PROOF. The required inequality

$$\inf_E \max(P(E), Q(\Omega \setminus E))$$
$$\geq \frac{1 - \sup_E |P(E) - Q(E)|}{2}$$

follows from

$$\max(P(E), Q(\Omega \setminus E)) \geq \frac{1 - |P(E) - Q(E)|}{2} \quad \forall E,$$

and the last inequality is true even when the max is replaced by the arithmetic mean. □

It is easy to check that Assumption 4 is satisfied for sampling from the Gaussian family $N(\theta, \sigma^2)$ with known variance $\sigma^2$ and with $I(\theta) := \sigma^{-2}$. For this model,

$$\|P_{n,\theta_1} - P_{n,\theta_2}\| = 1 - 2\Phi\left(-\frac{|\theta_2 - \theta_1|\sqrt{n}}{2\sigma}\right),$$

again with an equality and without the need to subtract $\varepsilon$.

## 5. CONCLUSION

It is widely accepted that advances in computing have brought about deep changes in the theory and practice of statistics. However, the use of the theory of computing and, in particular, of its core notion of computability, has been very limited in the classical areas of statistics, such as parameter estimation and hypothesis testing. The notion of computability appears to be especially useful in questions of efficiency and superefficiency, where it allows us to delineate the class of statistical procedures that we would like to compete with. In particular, restricting ourselves to computable estimators, we can ask whether a given parameter point $\theta$ can be a point of superefficiency. Hotelling's and Hodges's examples show that, without this restriction, the answer is vacuous: any $\theta$ can be a point of superefficiency. With the restriction, the answer we gave in Section 3 is that superefficiency is impossible at noncomputable $\theta$, whereas "hyperefficiency" ($T_n = \theta$ for all $n$) is possible at computable $\theta$.

This paper only deals with the most classical aspects of superefficiency. It does not even touch multivariate regular statistical models, let alone models in which rates of convergence of the maximum likelihood estimates are different from $n^{-1/2}$ and nonparametric models. Can the notion of computability be usefully applied in these and other more complex cases? I hope the answer is positive.

## APPENDIX A: COMPUTABILITY

The first paper to propose a general notion of computability, and to claim that its notion of computability is general, was Church's [6] (1936). Church considered functions $F : \mathbb{N}^m \to \mathbb{N}$, where $\mathbb{N}$ is the set of all positive integer numbers; for now, we will be only interested in the case $m = 1$. On the one hand, he formally defined his class of computable ("effectively calculable," as he said) functions, and on the other hand, he put forward the informal thesis (often referred to as the Church thesis) that his formal notion is the formalization of our intuitive notion of computability. One of Church's arguments in favor of the Church thesis was that two natural but very different definitions of effectively calculable functions, Church and Kleene's $\lambda$-definability and Herbrand and Gödel's recursiveness, are equivalent.

The Church thesis was further boosted by Alan Turing's observation [46] that Church's effective calculability is equivalent to computability using a formal model of a computing device, nowadays known as the Turing machine. A similar computing device was introduced at the same time by Emil Post [35], and another, rather different one, was introduced later by Andrei Markov, Jr. [30]; both devices led to the same class of computable functions as the Turing machine.

In 1953 Andrei Kolmogorov [22], later joined by his student Vladimir Uspensky [23], carefully analyzed the notion of an algorithm and introduced its very general formalization. Kolmogorov and Uspensky's goal was to show that "the most general, for the current state of science, notion of an algorithm" (my translation) leads to the same class of computable functions. As they had expected, their formalization (along with several other definitions they considered but did not include in the paper) indeed turned out to be equivalent to the previous ones.

At this time, there is a consensus that the intuitive notion of computability for functions $F : \mathbb{N} \to \mathbb{N}$ is indeed captured by the numerous available equivalent definitions. This notion will be assumed to be known in the rest of this paper; precise definitions can be found in, for example, Rogers's classical book [41].

A set $A \subseteq \mathbb{N}$ is called *decidable* if the function

$$\chi(x) := \begin{cases} 1, & \text{if } x \in A, \\ 2, & \text{otherwise,} \end{cases}$$

is computable. A function $F : A \to B$, where $A$ and $B$ are decidable subsets of $\mathbb{N}$, is said to be *computable* if its extension $F' : \mathbb{N} \to \mathbb{N}$ defined as

$$F'(x) := \begin{cases} F(x), & \text{if } x \in A, \\ 1, & \text{otherwise,} \end{cases}$$

is computable.

Many familiar countable sets $X$, such as $\mathbb{N}^2$, the set $\mathbb{Q}$ of all rational numbers, the set of all open intervals $(a, b) \subseteq \mathbb{R}$ with rational end-points, etc., can be represented as "spaces of finite objects" (in the terminology of Shoenfield [43]) by fixing a canonical injection $\phi_X : X \to \mathbb{N}$ mapping $X$ onto a decidable subset of $\mathbb{N}$. For example, a popular bijection $\phi_{\mathbb{N}^2} : \mathbb{N}^2 \to \mathbb{N}$ is the Cantor pairing function; it turns $\mathbb{N}^2$ into a space of finite objects. The reader will be assumed to be familiar with such canonical injections $\phi_X$ for the standard spaces of finite objects $X$. Intuitively, $\phi_X(x)$ encodes $x \in X$ as a positive integer, and instead of working with finite objects $x \in X$ directly, we can work with their codes.

The computability of $F : X \to Y$, where $X$ and $Y$ are spaces of finite objects, is defined as the computability of $\phi_Y \circ F \circ \phi_X^{-1} : \phi_X(X) \to \phi_Y(Y)$. A set $A \subseteq X$, where $X$ is a space of finite objects, is said to be *recursively enumerable* if $A = F(\mathbb{N})$ for some computable function $F : \mathbb{N} \to X$.

## A.1 Computable Real Numbers

The main goal of Turing's paper [46] was, in fact, not the definition of computable functions but the definition of computable real numbers. Turing's definition was that a real number is computable if its decimal expansion is computable. There are many equivalent definitions. For example, a real number $t$ is computable if and only if there exists a computable function $F : \mathbb{N} \to \mathbb{N}$ such that $||t| - F(n)/n| \leq 1/n$ for all $n \in \mathbb{N}$. This notion of a computable real number is as uncontroversial as the notion of a computable function $F : \mathbb{N} \to \mathbb{N}$.

Theorem 1 talks about computability of two objects: the estimator $\{T_n\}$ and the parameter point $\theta$. We have just defined what the computability of $\theta$ means. The situation with $\{T_n\}$ is more complicated. Typically, $T_n : \mathbb{R}^n \to \mathbb{R}$, and the notion of computability of real-valued functions of real numbers is notoriously ill-defined. There is the "core" notion of a computably continuous function, to be discussed in Section A.3, but there is no consensus about the "right" definition for more general classes of functions. In the next subsection we define computable estimators in an *ad hoc* manner, in order to obtain a strong statement of Theorem 1.

## A.2 Computable Estimators

The theory of computability over the real numbers often uses "effective" (i.e., computable in some sense) versions of various topological notions, such as openness, closeness, continuity, etc. A set $A \subseteq \mathbb{R}$ is said to be *effectively open* if it is the union of a recursively enumerable set of open intervals with rational end-points. In other words, $A$ is effectively open if it can be represented in the form $A = \bigcup_i (a_i, b_i)$, where $(a_i, b_i)$, $i \in \mathbb{N}$, is a computable sequence of open intervals with rational end-points. Complements $\mathbb{R} \setminus A$ of effectively open sets $A \subseteq \mathbb{R}$ are called *effectively closed*. More generally, sets $A_x \subseteq \mathbb{R}$ indexed by $x \in X$, where $X$ is a space of finite objects, are said to be *effectively open uniformly in $x$* if they can be represented in the form $A_x = \bigcup_i (a_i(x), b_i(x))$, where $(i, x) \mapsto (a_i(x), b_i(x))$ is a computable function mapping $\mathbb{N} \times X$ to the set of open intervals with rational end-points. In this case the complementary sets $\mathbb{R} \setminus A_x$ are said to be *effectively closed uniformly in $x$*.

Now we can define our notion of a computable estimator. As in Sections 3 and 4, we consider a sequence of statistical models $\{P_{n,\theta}\}$, where $\theta$ ranges over an open interval $\Theta \subseteq \mathbb{R}$; the end-points of $\Theta$ are assumed computable (by definition, $-\infty$ and $\infty$ are computable). Let $\{T_n\}$ be an estimator for $\{P_{n,\theta}\}$. We say that $\{T_n\}$ is *computable* if, for each $\delta \in \mathbb{Q} \cap [0, 1)$, the closures

$$(22) \qquad \overline{\{q \in \Theta \mid P_{n,q}(|T_n - q| > u) \leq \delta\}}$$

are effectively closed uniformly in $n \in \mathbb{N}$ and $u \in \mathbb{Q} \cap [0, \infty)$. Intuitively, the inequality "$\leq$" in (22) means that $T_n$ is a good estimator when the true parameter point is $q$, with $\delta$ and $u$ determining how demanding our notion of "good" is. The closure of the set of such $q$ is required to be uniformly effectively closed. It seems obvious that this condition will be satisfied for estimators $\{T_n\}$ specified by an explicit procedure.

Notice that our definition of computability of an estimator $\{T_n\}$ is in fact a joint requirement on the estimator and $\{P_{n,\theta}\}$. Interestingly, it does not impose any computability restrictions on the sample spaces $\Omega_n$, which do not enter the definition explicitly.

It is easy to check that our definition of computability of $\{T_n\}$ agrees with the definition of computability of a parameter point $\theta \in \Theta$, in the sense that the two notions coincide when $T_n := \theta$ is a constant estimator. Indeed, in this case,

$$\overline{\{q \mid P_{n,q}(|T_n - q| > u) \leq \delta\}} = [\theta - u, \theta + u] \cap \Theta,$$

and the last family of closed sets are effectively closed uniformly in $u \in \mathbb{Q} \cap [0, \infty)$ if and only if $\theta$ is computable.

Let us now check that the proof of Theorem 1 goes through for our definition of computability of $\{T_n\}$. Without loss of generality, we assume that $c$ and $\varepsilon$ in the proof of Theorem 1 are rational numbers and that $a$ and $I$ are chosen in such a way that $\delta := a\Phi(-c\sqrt{I})$ [cf. (11)] is a rational number; we have already said that $(L, R)$ is an interval with rational end-points. The requirement (13) leaves us enough freedom to make $n$ a square (i.e., to make $n^{1/2}$ integer), assuming that the interval $(L, R)$ is sufficiently short. Therefore, the closure of the set of $q \in \Theta$ satisfying (14) is effectively closed uniformly in the squares $n \geq N$. Assuming that $(\theta_1, \theta_2)$ is an interval with rational end-points [this is true initially, for $(\theta_1, \theta_2) = (L, R)$], we can compute a new interval $(\theta_1', \theta_2') \supseteq \mathfrak{S}$ with rational end-points of length $|\theta_2' - \theta_1'| \leq (1 + \varepsilon)^{-1}|\theta_2 - \theta_1|$. (Just use the definition of effective closeness and the compactness of bounded closed intervals in $\mathbb{R}$.)

## A.3 Computable Continuity

This subsection discusses the traditional notion of computability over the reals going back to the work of Brouwer on the intuitionistic foundations of mathematics; see [31] for an excellent description. Grzegorczyk [15] showed that this traditional notion of computability is equivalent to several other definitions considered in literature. An advantage of his exposition is that it is firmly based on the standard foundations of mathematics. The term "computable continuity" (in the form "computable continuous") is Grzegorczyk's ([15], footnote on page 71), who believed that it is possible to introduce some kinds of computable real functions which are not continuous.

Intuitively, a function $F$ defined over the reals is computably continuous if we can compute $F(x)$ to an arbitrary accuracy when given $x$ to an arbitrary accuracy. This condition indeed implies the continuity of $F$: for example, the simplest discontinuous function

$$F(x) := \begin{cases} 1, & \text{if } x \geq 0, \\ 0, & \text{otherwise,} \end{cases}$$

can never be computed to accuracy $1/3$ at the point $x = 0$, no matter how accurately we know $x$. On the other hand, any explicitly given continuous function the reader is likely to come across will be computable.

We start from defining what it means for a sequence of statistical models $\{P_{n,\theta}\}$ to be computably continuous (in the topology of weak convergence). As before, we assume that $\theta$ ranges over an open interval $\Theta$ of the real line $\mathbb{R}$ with computable end-points, and, for concreteness, we also assume that $P_{n,\theta}$ is a probability measure on $\Omega_n = \mathbb{R}^n$.

A *basic set* in $\mathbb{R}^m$ is the product $\prod_{i=1}^{m}(a_i, b_i)$ of bounded open intervals with rational end-points. An *elementary set* in $\mathbb{R}^m$ is a finite union of basic sets. The family of basic sets and the family of elementary sets can be regarded as spaces of finite objects. A subset of $\mathbb{R}^m$ is *effectively open* if it is the union of a recursively enumerable family of basic sets. A function $F : \Theta \to \mathbb{R}$ is *computably lower semicontinuous* if the set $\{(\theta, t) \mid F(\theta) > t\}$ is effectively open. The uniform versions of effective openness and computable lower semicontinuity are defined as before. A sequence of statistical models $\{P_{n,\theta}\}$ is said to be *computably continuous* if the function $P_{n,\theta}(E)$ is computably lower semicontinuous in $\theta$ uniformly in $n \in \mathbb{N}$ and elementary sets $E \subseteq \Omega_n$. This is a weak condition; the statistical models usually found in statistics textbooks are computably continuous.

Fix a computably continuous sequence of statistical models $\{P_{n,\theta}\}$. Let $\{T_n\}$ be an estimator for $\{P_{n,\theta}\}$. It is *computably continuous* if both sets $\{(x,t) \in \Omega_n \times \mathbb{R} \mid T_n(x) > t\}$ and $\{(x,t) \in \Omega_n \times \mathbb{R} \mid T_n(x) < t\}$ are effectively open uniformly in $n \in \mathbb{N}$.

It is not difficult to check that all computably continuous estimators are computable in our sense [see (22)]. In fact, for computably continuous estimators the operation of closure in (22) is superfluous: already the sets $\{q \mid P_{n,q}(|T_n - q| > u) \leq \delta\}$ are effectively closed uniformly in $n$ and $u$.

### A.4 Computability with an Oracle

The important idea of computability with an oracle was introduced by Turing [47]. An oracle Turing machine is allowed to read a tape containing an infinite sequence $S$ of symbols, not necessarily computable. Replacing in all our definitions computable functions $F : \mathbb{N} \to \mathbb{N}$ with $S$-computable functions $F : \mathbb{N} \to \mathbb{N}$ (i.e., functions computable by oracle Turing machines allowed to read $S$) leads to the notions of $S$-computable real numbers, $S$-computable estimators, $S$-computably continuous estimators, etc. Theorem 1 remains true if the two entries of "computable" are replaced by "$S$-computable." Since every estimator is $S$-computable for some $S$ (see Lemma 5 below), this "relativized" version of Theorem 1 contains Corollary 2 as a special case.

LEMMA 5. *Every estimator is $S$-computable for some $S$.*

PROOF. Let $(a_i, b_i)$, $i = 1, 2, \ldots$, be a computable enumeration of all open intervals with rational endpoints. It suffices to take as $S$ an infinite binary sequence encoding the function $F : \mathbb{N} \times \mathbb{N} \times (\mathbb{Q} \cap [0, \infty)) \times (\mathbb{Q} \cap [0, 1)) \to \{0, 1\}$ defined by the requirement that $F(i, n, u, \delta) = 1$ if and only if $(a_i, b_i)$ and (22) are disjoint. □

## APPENDIX B: DIRECT PROOF OF COROLLARY 2

This appendix gives a proof of Corollary 2 that does not use the notion of computability. It parallels the proof of Theorem 1.

Inequality (8) implies that there exist $a \in (0, 1)$, $I > I(\theta)$, and $N$ such that inequality (11) holds for all $n \geq N$. Since $a$ and $I$ can be taken rational, it suffices to prove for fixed $a$, $I$ and $N$ that (11) holds for all $n \geq N$ only for countably many $\theta$.

Suppose $\theta$ satisfies (11) for all $n \geq N$. It suffices to prove that there exists an open interval $(L, R) \ni \theta$ such that $\theta$ is the only point in $(L, R)$ satisfying (11) for all $n \geq N$. Let $\varepsilon > 0$, $N$ and $O \ni \theta$ satisfy the condition in Lemma 1. Take any $(L, R) \subseteq O$ satisfying

$$
(23) \qquad \frac{4(1+\varepsilon)^4 c^2}{(R-L)^2} \geq N,
$$

$$
\frac{4(1+\varepsilon)^6 c^2}{(R-L)^2} - \frac{4(1+\varepsilon)^4 c^2}{(R-L)^2} > 1,
$$

and $I(q) < I$ for all $q \in (L, R)$.

Suppose $(L, R)$ contains two distinct points $\theta_1$ and $\theta_2$ satisfying (11) for all $n \geq N$. Choose positive integer $n$ satisfying (12) and $n \geq N$. Such an $n$ exists since (12) is equivalent to (13) and we have assumed (23). By Lemma 1, (12) and (14) (applied to $q = \theta_1$ and $q = \theta_2$), there exists $\omega \in \Omega_n$ such that

$$
|T_n(\omega) - \theta_1| \leq cn^{-1/2},
$$
$$
|T_n(\omega) - \theta_2| \leq cn^{-1/2};
$$

therefore, the triangle inequality and (12) imply

$$
|\theta_2 - \theta_1| \leq 2cn^{-1/2} \leq (1+\varepsilon)^{-2}|\theta_2 - \theta_1|,
$$

which is impossible.

## ACKNOWLEDGMENTS

I am grateful to John Aldrich and Stephen Stigler for helpful comments. The paper very much benefitted from close reading by two referees and an Associate Editor. This work was supported in part by EPSRC (Grant EP/F002998/1).